\input amstex
\documentstyle{amsppt}
\loadbold
\topmatter
\title  Geometric Realization and  \\
 K-Theoretic Decomposition of C*-algebras
      \endtitle
\rightheadtext{  Geometric Realization of C*-algebras    }
\author    C. L. Schochet       \endauthor

\medbreak
\affil
Mathematics Department \\
Wayne State University \\
Detroit, MI 48202
\endaffil
\email     {    claude\@math.wayne.edu }     \endemail
\date {Submitted November 25, 1999 } \enddate
 \keywords     {$K$-theory for $C^*$-algebras, geometric realization,
Kasparov theory}
\endkeywords
\subjclass {  Primary 46L80, 19K35, 46L85}
\endsubjclass
\abstract {
Suppose that $A$ is a separable $C^*$-algebra and that $G_*$ is
a (graded) subgroup of the ${\Bbb Z} /2$-graded group $K_*(A)$. Then there
is a natural short exact sequence
$$
0 \to G_* \longrightarrow  K_*(A)  \longrightarrow K_*(A)/G_*  \to 0.
\tag *
$$
In this note we demonstrate how to geometrically realize this
sequence at the level of $C^*$-algebras. As a result,  we
$KK$-theoretically decompose $A$ as
$$
0 \to   A\otimes \Cal K
\longrightarrow A_f
\longrightarrow SA_t
\to 0
$$
where $K_*(A_t)$ is the torsion subgroup of $K_*(A)$ and
$K_*(A_f)$ is its torsionfree quotient. Then we further
decompose $A_t$: it is $KK$-equivalent to $\oplus _p A_p$
where $K_*(A_p)$ is the $p$-primary subgroup of the
torsion subgroup of $K_*(A)$.
 We then apply this
realization to study the Kasparov group $K^*(A)$ and related objects.

  }\endabstract
\endtopmatter

\document
\magnification = 1200
\pageno=1

\def\tensor{\otimes}

\def\KK #1.#2 { KK_*( #1,#2) }
 \def\KKgraded #1.#2.#3 { KK_{#1}( #2, #3) }

  \def\invlim{\underset \longleftarrow\to {lim} \, }

\def\ext #1.#2 {Ext _{\Bbb Z}^1( {#1} , {#2} ) }
\def\pext #1.#2 {Pext _{\Bbb Z}^1( {#1} , {#2} ) }
\def\pext #1.#2 {Pext _{\Bbb Z}^1( {#1} , {#2} ) }
\def\hom #1.#2 {Hom_{\Bbb Z} (#1 ,  #2  ) }

\def\dd{\bold X}

In Section 1 we produce the basic
geometric realization.
 For any separable
$C^*$-algebra $A$  and group $G_*$
we produce associated $C^*$-algebras $A_s$  ($s$ for subgroup)
and $A_q$ ($q$ for quotient group)
 and, most importantly, a short exact sequence of $C^*$-algebras
$$
0 \to A\tensor \Cal K  \to A_q  \to SA_s  \to 0
$$
whose associated $K_*$-long exact sequence is (*).
In the case where $G_*$ is the torsion subgroup of $K_*(A)$ we
use the notation $A_t$ ($t$ for torsion) and $A_f$ ($f$ for torsionfree)
 respectively. We further decompose $A_t$ into its $p$-primary
summands $A_p$ for each prime $p$.

Section 2 deals with the following question: may calculations
of the Kasparov groups $\KKgraded *.A.B $ be reduced down to
the four cases $(A_t,B_t)$, $(A_t, B_f)$, $(A_f, B_t)$ and
$(A_f, B_f)$ ? We show that this is indeed possible in a wide
variety of situations.
Sections 3 and 4 deal with these special cases.

Geometric realization as a general technique was introduced to
topological $K$-theory of spaces
by M. F. Atiyah
\cite{1}
in his proof of the K\" unneth theorem for $K^*(X \times Y)$. We adapted
the technique \cite{6}
 to prove the corresponding theorem for the $K$-theory for
$C^*$-algebras and used it with J. Rosenberg in our proof of the
Universal Coefficient Theorem (UCT) \cite{4}.
\medbreak

\beginsection {1. Geometric Realization}

In this section we produce the main geometric realization and we
extend the result to give a $p$-primary decomposition for a $C^*$-algebra.

\medbreak
Let $\Cal N$ denote the bootstrap category \cite {6, 4}.

\medbreak

\proclaim {Theorem 1.1}  Suppose that $A$ is a
separable $C^*$-algebra.
Let $G_*$ be some subgroup of $K_*(A)$.
Then there is an
associated $C^*$-algebra $A_s \in \Cal N$, a
separable  $C^*$-algebra $A_q$, and a short exact
sequence
$$
0 \to A\tensor \Cal K  \to A_q  \to SA_s \to 0
\tag 1.2
$$
whose induced $K$-theory long exact sequence fits
into the commuting diagram
$$
\CD
0 @>>> K_*(A_s)
@>>> K_*(A\tensor \Cal K) @>>> K_*(A_q) @>>> 0 \\
@.   @VV\cong V   @VV\cong V   @VV\cong V   \\
0 @>>> G_*
@>>> K_*(A) @>>> K_*(A)/G_*  @>>> 0.
\endCD
\tag 1.3
$$
If $A$ is nuclear then so is $A_q$. If $A \in \Cal N$ then so is $A_q$.
If $A \in \Cal N$ and if $G_*$ is a direct summand of $K_*(A)$ then
$A$ is $KK$-equivalent to $A_s \oplus A_q $.
\endproclaim

\medbreak
Note that we think of $A_s$ as realizing the subgroup $G_*$
and $A_q$ as realizing the quotient group $K_*(A)/G_*$,
hence the notation.

\demo{Proof} Let $A_s$ denote any $C^*$-algebra in $\Cal N$
 with
$$
K_*(A_s) \cong G_* .
$$
Such $C^*$-algebras exist and are unique up to $KK$-equivalence
 by the UCT \cite {4}. Let
$$
\theta : K_*(A_s) \to K_*(A)
$$
be the corresponding homomorphism. Since $A_s \in \Cal N$, the UCT
 holds for the pair $(A_s, A)$, and so
 $\theta $ is in the image of the index map
 $$
 \gamma : \KKgraded *.{A_s}.A   \to \hom {K_*(A_s)}.{K_*(A)}   .
 $$
Say that
$$
\theta  = \gamma (\tau )
$$
for some
$$
\tau \in \KKgraded 0.{A_s}.A      .
$$
As $A_s$ is nuclear,
$$
\KKgraded 0.{A_s}.A  \cong \Cal Ext(SA_s,A)
$$
and hence $\tau $ corresponds to an equivalence class of extensions of $C^*$-algebras of the form
$$
0 \to A\tensor \Cal K   \to  E  \to SA_s  \to 0.
$$
Define $A_q = E$. (This choice depends upon the choice of $A_s$
 among its $KK$-equivalence class and the choice of $\tau $ modulo the
  kernel of $\gamma $ ).
Note that $E$ is nuclear/bootstrap if and only if $A$ is nuclear/bootstrap.
  Then the diagram
$$
\CD
@.
K_j(A_q) @>>>
K_j(SA_s)  @>\delta >>   K_{j-1}(A\tensor \Cal K) @>>> K_{j-1}(A_q)  \\
@.@.
   @VV\cong V     @VV\cong V   \\
@.@. K_{j-1}(A_s)  @>\theta >>  K_{j-1}(A)
\endCD
$$
commutes, and thus $\delta $ is mono and
 the long exact  $K_*$-sequence breaks apart as shown.

If $G_*$ is a direct summand of $K_*(A)$ then
$$
K_*(A)  \cong G_*    \oplus K_*(A)/G_*
 \cong K_*(A_s) \oplus K_*(A_q) \cong K_*(A_s \oplus A_q)
$$
and,  replacing algebras by
 their suspensions as needed, the $KK$-equivalence
is obtained.  \enddemo\qed
\medbreak

Henceforth we shall regard $A_s$ and $A_q$ as $C^*$-algebras
 associated to $A$ and $G_*$, with the understanding that these
  are well-defined only up to $KK$-equivalence modulo the
kernel of $\gamma $, as explained above.

The next step is to decompose $A_t$ into its $p$-primary components.

\medbreak

\proclaim{Theorem 1.4} Let $A\in \Cal N$  and suppose
that $K_*(A)$ is a torsion group, so that $A = A_t$. Then $A$ is $KK$-equivalent to a $C^*$-algebra $\oplus A_p$, where
$$
K_*(A_p) \cong K_*(A)_p
$$
the $p$-primary torsion subgroup of $K_*(A)$.
\endproclaim

\medbreak

\demo{Proof} For each prime $p$, choose $N_{(p)} \in \Cal N$ with
$K_1(N_{(p)} ) = 0$ and
$$
K_0(N_{(p)} ) \cong \Bbb Z_{(p)}
$$
the integers localized at $p$. Define
$$
A_p = A_t \tensor N_{(p)}    .
$$
The K\"unneth formula \cite {6} implies that
$$
K_*(A_p) \cong K_*(A_t \tensor N_{(p)} ) \cong K_*(A_t) \tensor
K_*(N_{(p)} ) \cong K_*(A_t) \tensor \Bbb Z_{(p)}  \cong K_*(A)_p
$$
as desired. Then
$$
K_*(\oplus _p A_p) \cong \oplus _p  K_*(A_p)
\cong \oplus _p K_*(A)_p  \cong K_*(A_t)
$$
and another use of the UCT implies that $A_t$ is
$KK$-equivalent to $\oplus _p A_p$.
\enddemo\qed

\medbreak

We summarize:

\medbreak
\proclaim {Theorem 1.5}
Suppose that $A$ is a
separable $C^*$-algebra. Then there is an
associated $C^*$-algebra $A_t \in \Cal N$, a
separable  $C^*$-algebra $A_f$, and a short exact
sequence
$$
0 \to A\tensor \Cal K  \to A_f  \to SA_t \to 0
\tag 1.6
$$
whose induced $K$-theory long exact sequence fits
into the commuting diagram
$$
\CD
0 @>>> K_*(A_t)
@>>> K_*(A\tensor \Cal K) @>>> K_*(A_f) @>>> 0 \\
@.   @VV\cong V   @VV\cong V   @VV\cong V   \\
0 @>>> K_*(A)_t
@>>> K_*(A\tensor \Cal K) @>>> K_*(A)_f  @>>> 0.
\endCD
\tag 1.7
$$
If $A$ is nuclear then so is $A_f$. If $A \in \Cal N$ then so is $A_f$.
Further, the $C^*$-algebra $A_t $ has a $p$-primary decomposition:
 it  is $KK$-equivalent to a $C^*$-algebra $\oplus _p A_p$, where
 $A_p \in \Cal N$ for all $p$ and
$$
K_*(A_p) \cong K_*(A)_p
$$
the $p$-primary torsion subgroup of $K_*(A)$.
Finally, if $A \in \Cal N$ and $K_*(A)_t $ is a direct summand
of $K_*(A)$ then $A$ may be replaced  by
the $KK$-equivalent $C^*$-algebra
$A_t \oplus A_f$
\endproclaim\qed
\medbreak

\beginsection {2. Splitting the Kasparov Groups}

\medbreak

If $A$ and $B$ are in  $\Cal N$  and their $K$-theory
torsion subgroups
$K_*(A)_t$ and $K_*(B)_t$
 are direct summands then
 the final conclusion of
  Theorem 1.5 implies
that we may reduce the computation of $\KKgraded *.A.B $ to the
calculation of the four groups, namely
\roster
\item $\KKgraded *.{A_t}.{B_t}$
\medbreak
\item $\KKgraded *.{A_t}.{B_f}$
\medbreak
\item $\KKgraded *.{A_f}.{B_t}$
\medbreak
\item $\KKgraded *.{A_f}.{B_f}$  .
\endroster
\medbreak
We discuss the calculation of those groups in subsequent sections.
In this section we see what can be done {\it {without}}
assuming that the torsion subgroups are direct summands.

\medbreak
\proclaim {Theorem 2.1} Suppose that $A \in \Cal N$ and $K_*(B)$ is torsionfree. Then there is a short exact sequence
$$
0 \to  \KKgraded *.A_f .B  \to \KKgraded *.A.B \to \KKgraded *.A_t .B  \to 0     .
\tag 2.2
$$
In particular, letting $K^*(A) = \KKgraded *.A.{\Bbb C }$, there is a short exact sequence
$$
0 \to K^*(A_f)  \to K^*(A) \to K^*(A_t) \to 0    .
\tag 2.3
$$
If $K_*(B) $ is not necessarily torsionfree, then
 sequence 2.2 is exact if and only if the natural map
\medbreak
$$
\theta _h^* : \hom {K_*(A)}.{K_*(B)}  \to \hom {K_*(A_t)}.{K_*(B)}
\tag 2.4
$$
\medbreak\flushpar is onto, where $\theta : K_*(A_t) \to K_*(A) $
is the canonical inclusion.
\endproclaim
\medbreak\medbreak
Note that the map $\theta _h^*$ in (2.4) is frequently onto.
 This is the case, for instance, if $K_*(A_t)$ is a direct
  summand of $K_*(A)$.

The map $\theta $ is, up to isomorphism, the boundary
 homomorphism in the $K_*$-sequence associated to the short exact sequence
$$
0 \to A\tensor \Cal K  \to A_f  \to SA_t  \to 0
$$
and hence
$$
\theta (x) = x \tensor _{A_t}\delta
$$
where $\delta \in \KKgraded 1.{A_t}.A $ by \cite {9}.
Thus the map $\theta _h^*$ of (2.4) is induced from a $KK$-pairing.

\medbreak
\demo{Proof} Consider the commuting diagram
\medbreak
$$
\minCDarrowwidth {.2in}
\CD
@.   \hom{K_*(A_t)}.{K_*(B)}     @.          \KKgraded {*-1}.{A_t}.B      @.     0   \\
@.  @V{\beta}VV          @VVV         @VVV   \\
0 @>>>  \ext{K_*(A_f)}.{K_*(B)}  @>>>  \KKgraded *.{A_f}.B   @>>>   \hom{K_*(A_f)}.{K_*(B)}   @>>>   0  \\
@.  @VVV   @VVV   @VVV   \\
0 @>>>  \ext{K_*(A)}.{K_*(B)}  @>>>  \KKgraded *.{A}.B   @>>>   \hom{K_*(A)}.{K_*(B)}   @>>>   0  \\
@.  @V{\theta _e^*}VV      @V{\theta ^*}VV    @V{\theta _h^*}VV    \\
    0 @>>>  \ext{K_*(A_t)}.{K_*(B)}  @>>>  \KKgraded *.{A_t}.B
 @>>>   \hom{K_*(A_t)}.{K_*(B)}   @>>>  0   \\
@.   @VVV     @VVV   @V{\beta}VV \\
@.   0   @.   \KKgraded {*+1}.{A_f}.B     @.   \ext{K_*(A_f)}.{K_*(B)}
\endCD
$$
\vglue .3in
\flushpar  The three middle rows are exact by the UCT,
the middle column is exact by the exactness of $KK$, and the
two outer columns are exact by the standard $Hom$-$Ext$-sequence.

Suppose that $K_*(B)$ is torsionfree. Then
$$
\hom{K_*(A_t)}.{K_*(B)} = 0
\tag 2.5
$$
 since $ K_*(A_t) $ is a torsion group, and the surjectivity of
 $\theta _e^*$ implies the surjectivity of $\theta ^*$.

 If $K_*(B)$ is not necessarily torsionfree, then
the Snake Lemma \cite {11} implies that there is an exact sequence
$$
0 = Coker(\theta _e^*)  \longrightarrow Coker(\theta ^*)
\longrightarrow Coker(\theta _h^* ) \to 0
$$
and hence $\theta _h^*$ is onto if and only if $\theta ^*$ is
onto. The theorem then follows immediately, for the middle column
of the diagram degenerates to (2.2) if and only if $\theta ^*$ is
onto. \qed\enddemo

   Next we consider the dual situation, when $K_*(A)$ is a torsion
   group.

 \proclaim {Theorem 2.6} Suppose that $A \in \Cal N$ and that $K_*(A)$
  is a torsion  group. Then there is a natural exact sequence
  $$
  0 \to \KKgraded *.A.{B_t} \to \KKgraded *.A.B  \to \KKgraded *.A.{B_f}
   \to 0.
   \tag 2.7
  $$
  \medbreak\flushpar
  If $K_*(A)$ is not a torsion group then sequence
  (*) is exact if and only if the natural map
  $$
  \pi _* : \hom {K_*(A)}.{K_*(B)}  \to \hom {K_*(A)}.{K_*(B_f)}
  $$
  is onto, where $\pi : B\tensor \Cal K \to B_f $ is the natural map.
  \endproclaim

  \medbreak
  The proof of this result is dual to that of Theorem 2.1 and is omitted
for brevity.  \qed

\medbreak

\beginsection {3. Computing $\KKgraded *.{A_f}.B      $  }

In this section we consider the case where $K_*(A)$ is torsionfree
 (so that $A = A_f$). Recall \cite {2, 12} that
 a subgroup $H$ of an abelian group $K$ is {\it {pure}} if for
each positive integer $n$,
$$
nH = H \cap nG ,
$$
and an extension of groups
$$
0 \to H \to K \to G \to 0
$$
 is {\it {pure}} if $H$ is a pure subgroup of $K$.
 For abelian groups $G$ and $H$, $\pext G.H $ is
 the subgroup of $\ext G.H $ consisting of pure extensions.

\medbreak

Recall \cite {5, 8} that there is a natural topology
on the Kasparov groups and that with respect to this topology
the UCT sequence splittings constructed in \cite {4} are
continuous, so that the splitting is a splitting of topological
groups \cite{9}.
\medbreak

\proclaim {Theorem 3.1} Suppose that $A \in \Cal N$ and
 that $K_*(A)$ is torsionfree. Then there is a natural
  sequence of topological groups
$$
0 \to \pext {K_*(A)}.{K_*(B)} \to \KKgraded *.A.B
\to \hom {K_*(A)}.{K_*(B)}   \to 0
$$
The group
$
\pext {K_*(A)}.{K_*(B)}
 $
  is the closure of zero in the
natural topology on the group $\KKgraded *.A.B $ and
thus  the group
\newline
$\hom {K_*(A)}.{K_*(B)}$    is the Hausdorff quotient of
$\KKgraded *.A.B $.
\endproclaim

\medbreak

\demo{Proof} The UCT gives us the sequence
$$
0 \to \ext {K_*(A)}.{K_*(B)} \to \KKgraded *.A.B
\to \hom {K_*(A)}.{K_*(B)}   \to 0
$$
which splits unnaturally. If $K_*(A)$ is torsionfree then
$$
 \pext {K_*(A)}.{K_*(B)}
\,\cong\,    \ext {K_*(A)}.{K_*(B)}   .
$$

The remaining part of the theorem holds since
 we have shown in general \cite {10} that the group
$ \pext {K_*(A)}.{K_*(B)} $
 is the closure of
 zero in the natural topology on $\KKgraded *.A.B $ in
  the presence of the UCT. \enddemo\qed

\medbreak

We note that the resulting algebraic problems are
frequently very difficult.
 If $G$ is a torsionfree abelian group then $\hom G.H $
is unknown in general, though there is much known in special cases
(cf. \cite {2, 3}). The group $\pext G.H $ is also difficult,
though the case $\pext G.{\Bbb Z} $ is known (cf. \cite {3}). We
discuss $Pext$ in detail in \cite {12}.

\beginsection {4.  Computing $\KKgraded *.{A_t}.B $          }

\medbreak

In this section we concentrate upon the situation
 when $K_*(A)$ is a torsion group. Before beginning,
  we digress slightly to recall \cite {7} in more detail how one
   introduces coefficients into $K$-theory.

Given a countable abelian group  $G$,
 select some $C^*$-algebra $N_G \in \Cal N$ with
$$
K_0(N_G) = G \qquad\qquad K_1(N_G) = 0  .
$$
The $C^*$-algebra $N_G$ is unique up to $KK$-equivalence,
 by the UCT. Then for any $C^*$-algebra $A$,  define
$$
K_j(A; G) \cong  K_j (A\tensor N_G) .
\tag 4.1
$$

The K\"unneth Theorem \cite {6}
 implies that there is a natural short exact sequence
$$
0 \to K_j(A)\tensor G  \overset\alpha\to\longrightarrow
 K_j(A ; G) \to Tor_1^{\Bbb Z} (K_{j-1}(A) , G ) \to 0
\tag 4.2
$$
which splits unnaturally. If $G$ is torsionfree then $\alpha $ is an isomorphism
$$
\alpha :  K_j(A)\tensor G  \overset\cong\to\longrightarrow  K_j(A ; G).
$$

Let $\dd (G) = Hom (G, {\Bbb R}/{\Bbb Z}) $
 denote the Pontryagin dual of the group $G$.

\medbreak
\proclaim {Theorem 4.3} Suppose that $A \in \Cal N$ with $K_*(A)$
 a torsion group
  and suppose that $K_*(B)$ is torsionfree,
  so that $A = A_t$ and $B = B_f$.
 Then:
\medbreak
\roster
\item
$$
\KKgraded *.A.B \cong \ext {K_*(A)}.{K_{*-1}(B)}   .
$$
\medbreak
\item
$$
\KKgraded *.A.B \cong \hom {K_*(A)}.{K_{*-1}(B)\tensor \Bbb Q/\Bbb
Z }   .
$$
\medbreak
\item The group $\KKgraded *.A.B  $ is reduced and algebraically compact.
\medbreak
\item $$K^j(A) \cong \dd (K_{j-1}(A))  $$.
\medbreak
\item More generally, if $K_*(B)$ is finitely generated free, then
$$
\KKgraded j.A.B \cong \oplus _n \dd (K_{j-1}(A) )
$$
where $n$ is the number of generators of $K_*(B)$.
\endroster
\medbreak
\endproclaim

\medbreak

\demo {Proof} Part 1) follows at once from the UCT and the fact that there are no non-trivial homomorphisms from a torsion group to a torsionfree group. Part 2) follows from Part 1) by elementary homological algebra. Part 3) follows easily from a deep result of Fuchs and Harrison \cite {cf. 2, 46.1}:
 if $G$ is a torsion group then any group of the form $\hom G.H $ is reduced and algebraically compact. Part 4) follows from part 3) by setting $B = \Bbb C$ and observing that for any torsion  group $G$, we have
$$
\dd (G) = \hom G.{\Bbb Q / \Bbb Z }.
$$
\qed\enddemo

\medbreak

There is one additional case that fits into the present discussion and which partially overlaps with the result above.

\medbreak
\proclaim {Theorem  4.4} Suppose that $A \in \Cal N $
 and that $K_*(A) $ has no free direct summand.
  Then there is a natural short exact sequence of topological groups
$$
0  \to  \hom {K_*(A)}.{\Bbb R}        \to \dd (K_*(A))
\overset\chi\to\longrightarrow  K^*(A)  \to         0      .
\tag 4.5
$$
The map  $\chi : \dd (K_*(A))
\to  K^*(A)  $
is a degree one continuous open surjection. It is a homeomorphism if and only if $K_*(A)$ is a torsion group.
\endproclaim\

\medbreak

To be explicit about the grading,
 $$
\chi : \dd (K_j(A))
\to  K^{j-1}(A)
$$
which is the usual parity shift as torsion phenomena move from homology to cohomology.

\demo{Proof}  The UCT for $K^*(A) $ has the form
$$
0 \to \ext {K_*(A)}.{\Bbb Z}  \overset\delta\to\rightarrow  K^*(A)
\to \hom {K_*(A)}.{\Bbb Z}  \to  0
$$
with $\delta $ of degree one, so it suffices to compute $Ext$.
 In general the short exact sequence
$$
0 \to \Bbb Z \to \Bbb R  \to \Bbb R / \Bbb Z   \to 0
$$
yields a long exact sequence
$$
\hom {K_*(A)}.{\Bbb Z} \to \hom {K_*(A)}.{\Bbb R} \to
\dd (K_*(A))  \to \ext {K_*(A)}.{\Bbb Z}   \to 0.
$$
The fact that $K_*(A)$ has no free direct summand implies that  $\hom {K_*(A)}.{\Bbb Z} = 0$, so the
sequence degenerates to
$$
0 \to \hom {K_*(A)}.{\Bbb R} \to
\dd (K_*(A))  \to \ext {K_*(A)}.{\Bbb Z}   \to 0.
$$
Applying the UCT one obtains the sequence 4.5 as desired.
 The map $\chi $ is the composite of the UCT map and a natural
  homeomorphism. The rest of the Theorem is immediate.\enddemo\qed

\end

dvitps -d forward prim111999b.dvi |lpr -Pleonardo

\end